# ESTIMATING THE PROPORTION OF FALSE NULL HYPOTHESES AMONG A LARGE NUMBER OF INDEPENDENTLY TESTED HYPOTHESES

By Nicolai Meinshausen and John Rice

*ETH Zürich and University of California, Berkeley*

We consider the problem of estimating the number of false null hypotheses among a very large number of independently tested hypotheses, focusing on the situation in which the proportion of false null hypotheses is very small. We propose a family of methods for establishing lower $100(1-\alpha)\%$ confidence bounds for this proportion, based on the empirical distribution of the $p$-values of the tests. Methods in this family are then compared in terms of ability to consistently estimate the proportion by letting $\alpha \to 0$ as the number of hypothesis tests increases and the proportion decreases. This work is motivated by a signal detection problem that occurs in astronomy.

**1. Introduction.** An example that motivated our work is afforded by the Taiwanese–American Occultation Survey (TAOS), which we now briefly describe. The TAOS will attempt to detect small objects in the Kuiper Belt, a region of the solar system beyond the orbit of Neptune. The Kuiper Belt contains an unknown number of objects (KBOs), most of which are believed to be so small that they do not reflect enough light back to Earth to be directly observed. The purpose of the TAOS project is to estimate the number of these KBOs down to the typical size of cometary nuclei (a few kilometers) by observing occultations. The idea of the occultation technique is simple to describe. One monitors the light from a collection of stars that have angular sizes smaller than the expected angular sizes of comets. An occultation is manifested by detecting the partial or total reduction in the flux from one of the stars for a brief interval when an object in the Kuiper Belt passes between it and the observer. Four dedicated robotic telescopes will automatically monitor 2000–3000 stars every clear night for several years and their combined results will be used to test for an occultation of each star approximately every 0.20 seconds, yielding on the order of $10^{11}$ tests per year. The number of occultations expected per year ranges from tens to a few thousands, depending on what model of the Kuiper Belt is used. Having conducted a large number of tests, it is then of interest to estimate the number of occultations, or the occultation rate, since this will provide information on the distribution of KBOs. Note that in this context we are not so much interested in which particular null hypotheses are false as in how many are. The TAOS project was further described by Liang et al. [8] and Chen et al. [3].

We will base our analysis on the distribution of the $p$-values of the hypothesis tests. Let $\{G_\theta, \theta \in \Theta\}$ be some family of distributions, where $\theta$ is possibly infinite-dimensional and $G_0(t) = t$ with $0 \in \Theta$ is the uniform distribution on $[0, 1]$. All $p$-values are assumed to be independently distributed according to

$$P_i \sim G_{\theta_i}, \qquad i = 1, \ldots, n.$$

If a null hypothesis is true, the distribution of its $p$-value is uniform on $[0, 1]$ and $P_i \sim G_0$. We suppose that neither the family $\{G_\theta(t), \theta \in \Theta\}$ nor the parameter vector $(\theta_1, \ldots, \theta_n)$ is known, except from the fact that $G_0$ corresponds to the uniform distribution.

The proportion of null hypotheses that are false (the fraction of occultations in the TAOS example) is denoted by

$$\lambda = n^{-1} \sum_{i=1}^{n} \mathbf{1}\{\theta_i \neq 0\}. \tag{1}$$

Our goal is to construct a lower bound $\hat{\lambda}$ with the property

$$P(\hat{\lambda} \leq \lambda) \geq 1 - \alpha \tag{2}$$

for a specified confidence level $1 - \alpha$. Such a lower bound would allow one to assert, with a specified level of confidence, that the proportion of false null hypotheses is at least $\hat{\lambda}$. The global null hypothesis that there are no false null hypotheses can be tested at level $\alpha$ by rejecting when $\hat{\lambda} > 0$.

Our construction is closely related to that by Meinshausen and Bühlmann [9], which treats the case of possibly dependent tests, but with an observational structure that allows the use of permutation arguments that are not available in our case. Another estimate was examined by Nettleton and Hwang [10], but it does not have a property like (2). Our methodology is related to that of controlling the false discovery rate [1, 13], but the goals are different—we are not so much interested in which particular hypotheses are false as in how many are. However, we note that an estimate of the number of the false null hypotheses can be usefully employed in adaptive control of the false discovery rate [2]. In a modification of the original FDR method, Storey [13] also estimated the proportion of false hypotheses. The empirical distribution of $p$-values was used by Schweder and Spjøtvoll [11] to estimate the number of true null hypotheses; the methods used there are different than ours and do not provide explicit lower confidence bounds. The methods in this paper extend a proposal of Genovese and Wasserman [7]. We also relate our results to those of Donoho and Jin [6].

**2. Theory and methodology.** The estimate hinges on the definition of bounding functions and bounding sequences.

Let $U$ be uniform on $[0, 1]$. Let $U_n(t)$ be the empirical cumulative distribution function of $n$ independent realizations of a random variable with distribution $U$. For any real-valued function $\delta(t)$ on $[0, 1]$ which is strictly positive on $(0, 1)$, define $V_{n,\delta}$ as the supremum of the weighted empirical distribution

$$(3) \qquad V_{n,\delta} := \sup_{t \in (0,1)} \frac{U_n(t) - t}{\delta(t)}.$$

DEFINITION 1. A bounding function $\delta(t)$ is any real-valued function on $[0, 1]$ that is strictly positive on $(0, 1)$. A series $\beta_{n,\alpha}$ is called a bounding sequence for a bounding function $\delta(t)$ if, for a constant level $\alpha$:

(a) $n\beta_{n,\alpha}$ is monotonically increasing with $n$;
(b) $P(V_{n,\delta} > \beta_{n,\alpha}) < \alpha$ for all $n$.

The definition of a bounding sequence depends neither on the unknown proportion of false null hypotheses nor on the unknown distribution $G(t)$ of $p$-values under the alternative.

One is interested in the case where a proportion $\lambda$ of all hypotheses are false null hypotheses. Denote the empirical distribution of $p$-values by

$$(4) \qquad F_n(t) := n^{-1} \sum_{i=1}^{n} \mathbf{1}\{P_i \leq t\}.$$

Estimating the proportion of false null hypotheses can be achieved by bounding the maximal contribution of true null hypotheses to the empirical distribution function of $p$-values. We give a brief motivation. Suppose for a moment that there are only true null hypotheses. The expected fraction of $p$-values less than or equal to some $t \in (0, 1)$ equals, in this scenario, $U(t) = t$. The realized fraction $U_n(t)$ is, on the other hand, frequently larger than $t$. However, using Definition 1, the probability that $U_n(t)$ is larger than $t + \beta_{n,\alpha}\delta(t)$ is bounded by $\alpha$ simultaneously for all values of $t \in (0, 1)$. The proportion of $p$-values in the given sample that are in excess of the bound $t + \beta_{n,\alpha}\delta(t)$ can thus be attributed to the existence of a corresponding proportion of false null hypotheses and $F_n(t) - t - \beta_{n,\alpha}\delta(t)$ is hence a low-biased estimate of $\lambda$. As the bound for the contribution of true null hypotheses holds simultaneously for all values of $t \in (0,1)$, a lower bound for $\lambda$ is obtained by taking the supremum of $F_n(t) - t - \beta_{n,\alpha}\delta(t)$ over the interval $(0,1)$. A refined analysis shows that an additional factor $1/(1-t)$ can be gained when estimating the proportion of false null hypotheses.

DEFINITION 2. Let $\beta_{n,\alpha}$ be a bounding sequence for $\delta(t)$ at level $\alpha$. An estimate for the proportion $\lambda$ of false null hypotheses is given by

(5) $$\hat{\lambda} = \sup_{t \in (0,1)} \frac{F_n(t) - t - \beta_{n,\alpha}\delta(t)}{1-t}.$$

This estimate is indeed a lower bound for $\lambda$, as shown in the following theorem.

THEOREM 1. Let $\beta_{n,\alpha}$ be a bounding sequence for $\delta(t)$ at level $\alpha$ and let $\hat{\lambda}$ defined by (5). Then

(6) $$P(\hat{\lambda} \leq \lambda) \geq 1 - \alpha.$$

PROOF. The distribution of $p$-values $F_n$ is bounded by $F_n(t) \leq \lambda + (1-\lambda)U_{n_0}(t)$, where $n_0 = (1-\lambda)n$ and $U_{n_0}(t)$ is the empirical distribution of $n_0$ independent Uniform$(0,1)$-distributed random variables. Thus

(7) $$P(\hat{\lambda} > \lambda) \leq P\bigg(\sup_{t \in (0,1)} \frac{\lambda + (1-\lambda)U_{n_0}(t) - t - \beta_{n,\alpha}\delta(t)}{1-t} > \lambda\bigg)$$

(8) $$= P\bigg(\sup_{t \in (0,1)} (1-\lambda)(U_{n_0}(t) - t) - \beta_{n,\alpha}\delta(t) > 0\bigg)$$

(9) $$= P\bigg(\sup_{t \in (0,1)} U_{n_0}(t) - t - \frac{n}{n_0}\beta_{n,\alpha}\delta(t) > 0\bigg).$$

Whereas $n\beta_{n,\alpha}$ is monotonically increasing, $n\beta_{n,\alpha}/n_0 \geq \beta_{n_0,\alpha}$ and the proof follows by property (b) in Definition 1. □

2.1. *Asymptotic control.* Instead of finite-sample control, it is sometimes more convenient to resort to asymptotic control. A sequence $\beta_{n,\alpha}$ is said to be an *asymptotic bounding sequence* if $\beta_{n,\alpha}$ satisfies condition (a) from Definition 1 and, additionally, a modified condition (b'),

(10) $$\limsup_{n \to \infty} P(V_{n,\delta} > \beta_{n,\alpha}) < \alpha,$$

where $V_{n,\delta}$ is defined as in (3). If we suppose that the absolute number of false null hypotheses $n\lambda$ is growing with $n$, that is, $n\lambda \to \infty$ for $n \to \infty$, then for an asymptotic bounding sequence,

$$\limsup_{n \to \infty} P(\hat{\lambda} \leq \lambda) \geq 1 - \alpha.$$

Asymptotic control is typically useful in the following situation. For a given bounding function $\delta(t)$ and two sequences $a_n, b_n$, consider weak convergence of

$$a_n V_{n,\delta} - b_n \xrightarrow{D} L \tag{11}$$

to a distribution $L$. Any sequence $\beta_{n,\alpha}$ that satisfies the monotonicity condition (a) of Definition 1 and, additionally, $\beta_{n,\alpha} \geq a_n^{-1}(L^{-1}(1-\alpha) + b_n)$, is thus an asymptotic bounding sequence at level $\alpha$.

As an important example, consider the bounding function $\delta(t) = \sqrt{t(1-t)}$. The following lemma is due to Jäschke and can be found in [12], page 599, Theorem 1 (18).

LEMMA 1. *Let $a_n = \sqrt{2n \log \log n}$ and $b_n = 2 \log \log n + \frac{1}{2} \log \log \log n - \frac{1}{2} \log 4\pi$. Then*

$$a_n \sup_{t \in (0,1)} \frac{U_n(t) - t}{\sqrt{t(1-t)}} - b_n \xrightarrow{D} E^2, \tag{12}$$

*where $E$ is the Gumbel distribution $E(x) = \exp(-\exp(-x))$.*

REMARK 1. The convergence in (12) is in general slow. Nevertheless, the result is of interest here. First, the number of tested hypotheses is potentially very large (e.g., $10^{12}$ in the TAOS setting described in the Introduction). Moreover, the slow convergence is mainly caused by values of $t$ that are of order $1/n$. The expected value of the smallest $p$-value of true null hypotheses is at least $1/n$ and it might be useful to truncate in practice the range over which the supremum is taken in (5) to $(1/n, 1-1/n)$. Doing so, the following asymptotic results are still valid, while the approximation by the Gumbel distribution is empirically a good fit even for moderate values of $n$ [6].

Similar weak convergence results for other bounding functions can be found in [4] or [12].

2.2. *Bounding functions.* The estimate is determined by the choice of the function $\delta(t)$, the so-called bounding function, and a suitable bounding sequence.

There are many conceivable bounding functions. Bounding functions of particular interest include:

– *linear bounding function* $\delta(t) = t$;
– *constant bounding function* $\delta(t) = 1$;
– *standard deviation–proportional bounding function* $\delta(t) = \sqrt{t(1-t)}$.

The *linear bounding function* is closely related to the false discovery rate (FDR), as introduced by Benjamini and Hochberg [1]. In the FDR setting, the empirical distribution of $p$-values is compared to the linear function $t/\alpha$. The last down-crossing of the empirical distribution over the line $t/\alpha$ determines the number of rejections that can be made when controlling FDR at level $\alpha$. It is interesting to compare this to the current setting. In particular, it follows by a result of Daniels [5] that

$$P\left(\sup_{t\in(0,1)} U_n(t)/t \geq \lambda\right) = 1/\lambda.$$

The optimal bounding sequence at level $\alpha$ is thus given for the linear bounding function by $\beta_{n,\alpha} = 1/\alpha - 1$. Let $\hat{\lambda}$ be the estimate under the linear bounding function. The estimate vanishes hence, that is, $\hat{\lambda} = 0$, if and only if no rejections can be made under FDR control at the same level. Note that the bounding sequence is independent of the number of observations. This leads to weak power to detect the full proportion $\lambda$ of false null hypotheses when the proportion $\lambda$ is rather high but the distribution of $p$-values under the alternative deviates only weakly from the uniform distribution, as shown in an asymptotic analysis below.

An estimate under a *constant bounding function* was already proposed by Genovese and Wasserman [7]. Using the Dvoretzky–Kiefer–Wolfowitz (DKW) inequality, a bounding sequence is given by $\beta_{n,\alpha}^2 = \frac{1}{2n} \log \frac{2}{\alpha}$. In contrast to the linear bounding function, this bounding function sequence vanishes for $n \to \infty$. However, the estimate is unable to detect any proportion of false null hypotheses that is of smaller order than $\sqrt{n}$. The intuitive reason is that the bounding function $\delta(t)$ is not vanishing for small values of $t$. Any evidence from false null hypotheses, however strong it may be, is hence lost if there are just a few false null hypotheses.

As already argued above, a bounding sequence for the *standard deviation–proportional bounding function* is given by

$$\beta_{n,\alpha} = a_n^{-1}(E^{-1}(1-\alpha) + b_n), \tag{13}$$

where $E$ is the Gumbel distribution and $a_n, b_n$ are defined as in Lemma 1. Note that the bounding sequence is vanishing at almost the same rate as for the constant bounding function. In contrast to the constant bounding function, however, the standard deviation–proportional bounding function vanishes for small $t$. It will be seen that the standard deviation–proportional bounding function possesses optimal properties among a large class of possible bounding functions.

2.3. *Asymptotic properties of bounding sequences.* Faced with an enormous number of potential bounding functions, it is of interest to look at general properties of bounding functions, especially the asymptotic behavior of the resulting estimates. The asymptotic properties turn out to be mainly determined by the behavior of $\delta(t)$ close to the origin.

DEFINITION 3. For every $\nu \in [0, 1]$, let $Q_\nu$ be a family of real-valued functions on $[0, 1]$. In particular, $\delta(t) \in Q_\nu$ iff:

(a) $\delta(t)$ is nonnegative and finite on $[0, 1]$ and strictly positive on $(0, 1)$;
(b) $\delta(1 - t) \geq \delta(t)$ for $t \in (0, \frac{1}{2})$;
(c) the function $\delta(t)$ is regularly varying with power $\nu$, that is,

$$\lim_{t \to 0} \frac{\delta(bt)}{\delta(t)} = b^\nu.$$

Most bounding functions of interest are members of $Q_\nu$ for some value of $\nu \in [0, 1]$. The constant bounding function is a member of $Q_0$, while the linear bounding function is a member of $Q_1$ and the standard deviation–proportional bounding function is a member of $Q_{1/2}$.

It holds in general for any bounding function that bounding sequences cannot be of smaller order than the inverse square root of $n$. In particular, note that by Definition 1 of a bounding sequence, it has to hold for any $t \in (0, 1)$ that $P(U_n(t) - t - \beta_{n,\alpha}\delta(t) > 0) < \alpha$ for all $n \in \mathbb{N}$. Whereas $nU_n(t) \sim B(n, t)$ is binomially distributed with mean $nt$ and variance proportional to $n$, it follows indeed that

$$\liminf_{n \to \infty} n^{1/2} \beta_{n,\alpha} > 0.$$

Consider now bounding functions $\delta(t)$, which are members of $Q_\nu$ with some $\nu \in (\frac{1}{2}, 1]$. It follows directly from Theorem 1.1(iii) in [4], page 255, that a more restrictive assumption has to hold in this case, namely

(14) $$\liminf_{n \to \infty} n^{1-\nu} \beta_{n,\alpha} > 0.$$

For $\nu = 1$ this amounts to $\liminf_{n \to \infty} \beta_{n,\alpha} > 0$. The linear bounding function is a member of $Q_1$, explaining the lack of convergence to zero of the corresponding optimal bounding sequence $1/\alpha - 1$.

For bounding functions $\delta(t) \in Q_\nu$ with $\nu \in [0, \frac{1}{2}]$, there exists some constant $c > 0$ so that $c\delta(t)^2 \geq t(1 - t)$. Hence, using Lemma 1, there exist bounding sequences so that

(15) $$\limsup_{n \to \infty} \left(\frac{n}{\log \log n}\right)^{1/2} \beta_{n,\alpha} < \infty.$$

The different asymptotic behavior of the bounding sequences influences the asymptotic power to detect false null hypotheses, as will be seen subsequently.

**3. Power.** We examine the influence of the bounding function $\delta(t)$ on the power to detect false null hypothesis. For simplicity of exposition, it is assumed that the $p$-values of all false null hypotheses follow a common distribution $G$, while $p$-values of true null hypotheses have a uniform distribution on $[0,1]$. For some $\gamma \in (0,1)$, let

$$\lambda \sim n^{-\gamma}.$$

A value of $\gamma = 0$ corresponds to a fixed proportion of false null hypotheses, while $\gamma = 1$ corresponds to a fixed absolute number of false null hypotheses. Here all cases between those two extremes are considered.

*Bounding sequences with vanishing level.* For the asymptotic analysis, it is convenient to let $\alpha = \alpha_n$ decrease monotonically for $n \to \infty$, so that $\alpha_n \to 0$ for $n \to \infty$. Note that $\alpha_n \to 0$ is equivalent to $P(V_{n,\delta} > \beta_{n,\alpha_n}) \to 0$ for $n \to \infty$. For notational simplicity, this assumption is strengthened slightly to

$$(16) \qquad V_{n,\delta}/\beta_{n,\alpha_n} \xrightarrow{p} 0, \qquad n \to \infty.$$

In almost all cases of interest, (16) and $\alpha_n \to 0$ are equivalent. To maintain reasonable power, one would like to avoid letting the level $\alpha_n$ vanish too fast as $n \to \infty$. For bounding functions $\delta(t) \in Q_\nu$ with $\nu \in [0, \frac{1}{2}]$ it is required that

$$(17) \qquad \limsup_{n \to \infty} \left(\frac{n}{\log n}\right)^{1/2} \beta_{n,\alpha_n} < \infty.$$

It follows from (15) that it is always possible to find a sequence $\alpha_n \to 0$ so that both (16) and (17) are satisfied. If both (16) and (17) are satisfied, the sequence $\alpha_n$ is said to *vanish slowly*. For bounding functions $\delta(t) \in Q_\nu$ with $\nu \in (1/2, 1]$, it will be seen below that the power is poor no matter how slowly the sequence $\alpha_n$ vanishes for $n \to 0$.

3.1. *Case* I: *many false null hypotheses,* $\gamma \in [0, \frac{1}{2})$. The fluctuations in the empirical distribution function are negligible compared to the signal from false null hypotheses if $\gamma \in [0, \frac{1}{2})$. Hence one should be able to detect (asymptotically) the full proportion of false null hypotheses in this first setting.

This is indeed achieved, as long as we look for bounding functions in $Q_\nu$ with $\nu \in [0, \frac{1}{2}]$, as shown below. If on the other hand $\nu \in (\frac{1}{2}, 1]$, one is in general unable to detect the full proportion of false null hypotheses. The proportion of detected false null hypotheses even converges in probability to zero for large values of $\gamma$ if $\nu$ is in the range $(\frac{1}{2}, 1]$. This includes in particular the linear FDR-style bounding function $t \in Q_1$, which is only able to detect a nonvanishing proportion of false null hypotheses (asymptotically) as long as the proportion $\lambda$ is bounded from below, which is only satisfied for $\gamma = 0$.

THEOREM 2. *Let $G$ be continuous and let $\inf_{t \in (0,1)} G'(t) = 0$. Let $\hat{\lambda}$ be the estimate under bounding function $\beta_{n,\alpha} \delta(t)$, where $\delta(t) \in Q_\nu$ with $\nu \in [0,1]$ and $\beta_{n,\alpha}$ is a bounding sequence. If $\nu \in [0, \frac{1}{2}]$ and $\alpha_n$ vanishes slowly, then, for all $\gamma \in [0, \frac{1}{2})$,*

$$\frac{\hat{\lambda}}{\lambda} \xrightarrow{p} 1, \qquad n \to \infty.$$

*However, for $\nu \in (\frac{1}{2}, 1]$ and $\gamma \in (1 - \nu, \frac{1}{2})$,*

$$\frac{\hat{\lambda}}{\lambda} \xrightarrow{p} 0, \qquad n \to \infty.$$

REMARK 2. The case $\inf_{t \in (0,1)} G'(t) = 0$ corresponds to the "pure" case in [7]. If $\inf_{t \in (0,1)} G'(t) > 0$, the results above (and below) hold if $\lambda$ is replaced by

$$\underline{\lambda} = \left(1 - \inf_{t \in (0,1)} G'(t)\right) \lambda.$$

Without making parametric assumptions about the distribution $G$ under the alternative, identifying $\underline{\lambda}$ is indeed the best one can hope for.

The message from Theorem 2 is that one should look for bounding functions in $Q_\nu$ with $\nu \in [0, \frac{1}{2}]$. This guarantees proper behavior of the estimate if the proportion $\lambda$ of false null hypotheses is vanishing more slowly than the square root of the number of observations.

3.2. *Case II: few false null hypotheses, $\gamma \in [\frac{1}{2}, 1)$.* As seen above, bounding functions in $Q_\nu$ with $\nu \leq \frac{1}{2}$ detect asymptotically the full proportion $\lambda$ of false null hypotheses if $\lambda$ is vanishing not as fast as the square root of the number of observations.

For $\gamma > 1/2$, no method can detect asymptotically the full proportion of false null hypotheses if the distribution under the alternative is fixed. For a fixed nondegenerate alternative, the majority of $p$-values from false null hypotheses fall with high probability into a fixed interval that is bounded away from zero. The fluctuations of the empirical distribution function in such an interval are asymptotically infinitely larger than any signal from false null hypotheses if $\gamma > 1/2$, which makes detection of the full proportion of false null hypotheses impossible.

It is hence interesting to consider cases where the signal from false null hypotheses is increasing in strength. Therefore, let $G = G^{(n)}$, the distribution of $p$-values under the alternative, be a function of the number $n$ of tests to conduct. The superscript is dropped in the following for notational simplicity.

*Shift-location testing.* It is perhaps helpful to think about $G$ as being induced by some shift-location testing problem. For each test it is assumed that there is a test statistic $Z_i$, which follows some distribution $T_0$ under the null hypothesis $H_{0,i}$ and some shifted distribution $T_{\mu_n}$ under the alternative $H_{1,i}$:

$$H_{0,i} : Z_i \sim T_0,$$
$$H_{1,i} : Z_i \sim T_{\mu_n}. \tag{18}$$

In the Gaussian case this amounts, for example, to $T_0 = \mathcal{N}(0,1)$ and $T_\mu = \mathcal{N}(\mu_n, 1)$. To have an interesting problem, one needs for $\gamma \in (\frac{1}{2}, 1)$ in general that the shift $\mu_n$ between the null and alternative hypotheses be increasing for an increasing number of tests; that is, $\mu_n \to \infty$ for $n \to \infty$.

On the other hand, one would like to keep the problem subtle. For the Gaussian case it was shown by Donoho and Jin [6] that an interesting scaling is given by $\mu_n = \sqrt{2r \log n}$ with $r \in (0,1)$. In this regime, the smallest $p$-value stems with high probability from a true null hypothesis. The false null hypotheses have hence little influence on the extremes of the distribution.

Instead of assuming Gaussianity of the test statistics, Donoho and Jin [6] considered a variety of different distributions. Under a generalized Gaussian (Subbotin) distribution, the density is for some positive value of $\kappa$ proportional to

$$T'_\mu(x) \propto \exp\left(-\frac{|x-\mu|^\kappa}{\kappa}\right).$$

The case $\kappa = 2$ corresponds clearly to a Gaussian distribution; $\kappa = 1$ corresponds to the double exponential case. The shift parameter is chosen then as

$$\mu_n = (\kappa r \log n)^{1/\kappa} \tag{19}$$

for some $r \in (0,1)$. Note that the expectation of the smallest $p$-value from true null hypotheses vanishes like $n^{-1}$, whereas under the scaling (19), the median $p$-value of false null hypotheses vanishes like $n^{-r}$ for $n \to \infty$ with some $r \in (0,1)$. In fact, consider for any member of the generalized Gaussian Subbotin distribution the $q$-quantile $G^{-1}(q)$ of the distribution of $p$-values under the alternative. For some constant $c_q$, the $q$-quantile is proportional to

$$G^{-1}(q) \propto \int_{\mu_n + c_q}^\infty \exp\left(-\frac{x^\kappa}{\kappa}\right) dx.$$

Applying l'Hôpital's rule twice, it follows for any $c$ and $\kappa > 0$ that

$$\lim_{a \to \infty} \frac{\log \int_{a+c}^\infty \exp(-x^\kappa/\kappa)\, dx}{-a^\kappa/\kappa} = 1.$$

Thus, under the scaling (19), for any every $q \in (0,1)$ and positive $\kappa$, the scaling of the $q$-quantile is given by

$$\log G^{-1}(q) \sim -r \log n. \tag{20}$$

With probability converging to 1 for $n \to \infty$, a $p$-value under a false null hypothesis is hence larger than the smallest $p$-value from all true null hypotheses as long as $r \in (0,1)$. For $r > 1$, the problem gets trivial as the probability that an arbitrarily high proportion of $p$-values under false null hypotheses is smaller than the smallest $p$-value from all true null hypotheses converges to 1 for $n \to \infty$.

The point of introducing the shift-location model under generalized Gaussian Subbotin distributions was just to identify (20) with $r \in (0,1)$ as the interesting scaling behavior of quantiles of $G$, the $p$-value distribution for alternative hypotheses. The setting (20) is potentially of interest beyond any shift-location model. We adopt the scaling (20) for the following discussion without making any explicit distributional assumptions about underlying test statistics.

THEOREM 3. *Let $\lambda \sim n^{-\gamma}$ with $\gamma \in [\frac{1}{2}, 1)$ and let the distribution $G$ of $p$-values under the alternative satisfy (20) for some $r \in (0,1)$. Let $\hat{\lambda}$ be the estimate of $\lambda$ under a bounding function $\beta_{n,\alpha}\delta(t)$, where $\delta(t) \in Q_\nu$ with $\nu \in [0, \frac{1}{2}]$ and $\beta_{n,\alpha_n}$ is a bounding sequence for $\delta(t)$. Let $\alpha_n$ vanish slowly. If $r > \frac{1}{\nu}(\gamma - \frac{1}{2})$,*

$$\frac{\hat{\lambda}}{\lambda} \xrightarrow{p} 1. \tag{21}$$

*If, on the other hand, $r < \frac{1}{\nu}(\gamma - \frac{1}{2})$, then*

$$\frac{\hat{\lambda}}{\lambda} \xrightarrow{p} 0. \tag{22}$$

REMARK 3. The analysis was only carried out for functions with $\nu \in [0, \frac{1}{2}]$ due to the deficits of the functions with $\nu \in (\frac{1}{2}, 1]$ discussed in the previous section. Nevertheless, it would be possible to carry out the same analysis here. For $\nu = 1$, one obtains, for example, a critical boundary $r > \gamma$.

The message from the last theorem is that among all bounding functions in $Q_\nu$ with $\nu \in [0, \frac{1}{2}]$, it is best to choose a member of $Q_{1/2}$. Bounding functions in $Q_{1/2}$ increase the chance to detect the full proportion $\lambda$ of false null hypotheses, as illustrated for a few special cases in Figure 1. The area in the $(r, \gamma)$ plane where $\hat{\lambda}/\lambda$ converges in probability to 1 for a bounding function in $Q_{1/2}$ includes in particular all areas of convergence for bounding functions in $Q_\nu$ with $\nu \in [0, \frac{1}{2}]$.

3.3. *Connection to the familywise error rate.* A different estimate of $\lambda$ is obtained by controlling the familywise error rate (FWER). In particular, let the estimate be the total number of $p$-values less than the FWER threshold $\alpha/n$, divided by the total number of hypotheses,

$$\hat{\lambda} = F_n\left(\frac{\alpha}{n}\right).$$

This is an estimate of $\lambda$ with the desired property $P(\hat{\lambda} > \lambda) < \alpha$. Controlling the familywise error rate has often been criticized for lack of power. Indeed, in the asymptotic analysis above it is straightforward to show that the area in the $(r, \gamma)$ plane where $\hat{\lambda}/\lambda \to_p 1$ is restricted to the half-plane $r > 1$ (neglecting again what happens directly on the border $r = 1$). In comparison to other estimates proposed here, the familywise error rate is hence particularly bad for estimating $\lambda$ if there are many false null hypotheses, each with a very weak signal. In addition, the construct requires that $p$-values can be determined accurately down to precision $\alpha/n$, which might be prohibitively small. In contrast, the performance of estimates of the form (5) does not deteriorate significantly if $p$-values are truncated at larger values.

The drawbacks of the familywise error rate are a consequence of the stricter inference one is trying to make when controlling the familywise error rate. In particular, one is trying to infer exactly *which* hypotheses are false nulls as opposed to only *how many* false nulls there are in total. The loss in power is hence the price one pays for this more ambitious goal.

3.4. *Connection to higher criticism.* A connection of the proposed estimate to the higher criticism method of Donoho and Jin [6] for detection of sparse heterogeneous mixtures emerges. In their setup $p$-values $P_i$, $i = 1, \ldots, n$, are i.i.d. according to a mixture distribution

$$P_i \sim (1 - \overline{\lambda})H + \overline{\lambda}G,$$

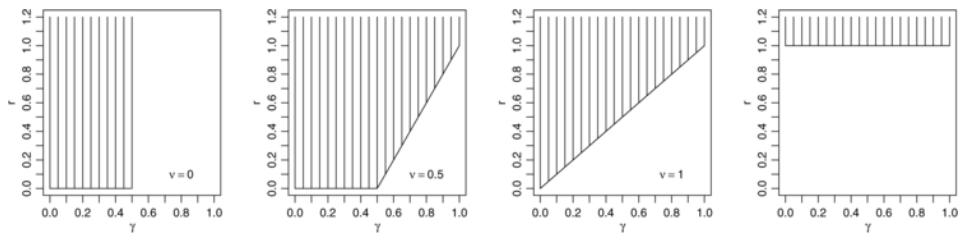

FIG. 1. *For $\nu = 0$ (left), $\nu = 1/2$ (second from left) and $\nu = 1$ (second from right), an illustration of the asymptotic properties of the estimate $\hat{\lambda}$. The shaded area marks those areas in the $(r, \gamma)$ plane where $\hat{\lambda}/\lambda \to_p 1$, whereas for the white areas $\hat{\lambda}/\lambda \to_p 0$. The choice $\nu = 1/2$ is seen to be optimal. The corresponding plot for control of the familywise error rate is shown on the right for comparison.*

where $H$ is the uniform distribution and $G$ the distribution of $p$-values under the alternative hypothesis. In [6] the focus is on testing the global null hypothesis that there are no false null hypotheses at all,

$$H_0: \overline{\lambda} = 0.$$

In contrast, in this current paper we are interested in quantifying the proportion $\lambda$ of false null hypotheses. The proportion $\lambda$ of false null hypotheses, as defined for the current paper in (1), can be viewed as a realization of a random variable with a binomial distribution, $n\lambda \sim B(n, \overline{\lambda})$. For the asymptotic considerations of this paper, however, the distinction between $\overline{\lambda}$ and $\lambda$ is of little importance because the ratio $\lambda/\overline{\lambda}$ converges almost surely to 1 for $n \to \infty$.

The two goals of higher criticism and the current paper are connected. If there is evidence for a positive proportion of false null hypotheses with the proposed method, then the global null $H_0$ can clearly be rejected. In other words, if one obtains a positive estimate $\hat{\lambda} > 0$ with $P(\hat{\lambda} > \lambda) < \alpha$, then the global null hypothesis $H_0: \overline{\lambda} = 0$ can be rejected at level $\alpha$. Note that the level is correct even for finite samples and not just asymptotically.

The connection between the two methods works as well in the reverse direction if an optimal bounding function is chosen. It emerged in particular from the analysis above that bounding functions that are members of $Q_{1/2}$ have optimal asymptotic properties. For the particular choice of a standard deviation–proportional bounding function in $Q_{1/2}$, let $\hat{\lambda}$ be an estimate of $\lambda$ and let $\beta_{n,\alpha}$ be a bounding sequence that satisfies

$$\beta_{n,\alpha} = n^{-1/2}(2\log\log n)^{1/2}(1 + o(1)).$$

Donoho and Jin [6] are not specific about choice of a critical value for higher criticism. However, choosing $\sqrt{n}\beta_{n,\alpha}$ as a critical value meets their requirements. The higher criticism procedure rejects in this case if and only if the estimate $\hat{\lambda}$ of the proportion of false null hypotheses is positive,

$$\{\text{reject } H_0: \overline{\lambda} = 0 \text{ with higher criticism}\} = \{\hat{\lambda} > 0\}.$$

If both $\overline{\lambda} \sim n^{-\gamma}$ and $\lambda \sim n^{-\gamma}$ for some $\gamma \in [0,1]$, the question arises if the area in the $(\gamma, r)$ plane where

(23) $$P(\text{higher criticism rejects } H_0) \to 1$$

is identical to the area where

(24) $$\frac{\hat{\lambda}}{\lambda} \xrightarrow{p} 1.$$

Intuitively, it is clear that it is somewhat easier to test for the global null hypothesis $H_0: \overline{\lambda} = 0$, as done in higher criticism, than to estimate the precise proportion $\lambda$ of false null hypotheses, as done in this paper. One would therefore expect that the area of convergence in the $(\gamma, r)$ plane of (23) includes the area of convergence of (24).

It is hence maybe surprising that for some cases the areas of convergence in the $(\gamma, r)$ plane of (23) and (24) agree. To illustrate the point, consider again the shift-location model (18) under a generalized Gaussian Subbotin distribution with parameter $\kappa \in (0, 2)$ and a shift (19) of test statistics under the alternative.

The area in the $(\gamma, r)$ plane where $\hat{\lambda}/\lambda \to_p 1$ is in this setting independent of the parameter $\kappa$. The detection boundary for higher criticism, however, does depend on $\kappa$. For the Gaussian case ($\kappa = 2$) and in general for $\kappa > 1$, the detection boundary for higher criticism is, for $\gamma \in (1/2, 1)$, below the area where $\hat{\lambda}/\lambda \to_p 1$. The reason for this is intuitively clear. The higher criticism method looks in these cases for evidence against $H_0$ in the extreme tails of the distribution $G$; see [6]. At these points, only a vanishing proportion of all $p$-values from false null hypotheses can be found. If one is trying to estimate the full proportion of false null hypotheses, the evidence for a certain amount of false null hypotheses has to be found at less extreme points, where one can expect a significant proportion of $p$-values from false null hypotheses. This limits the region of convergence in the sense of (24) compared to the area where higher criticism can successfully reject the global null hypothesis $H_0 : \overline{\lambda} = 0$.

However, for $\kappa \leq 1$ (including thus the case of a double-exponential distribution) the two areas where (23) and (24) hold, respectively, are identical, as shown in Figure 2. In the white area, both higher criticism and the current method fail to detect (asymptotically) the presence of false null hypotheses, and not even the likelihood ratio test is able to reject in these cases (asymptotically) the global null hypothesis $H_0 : \overline{\lambda} = 0$ that there are only true null hypotheses [6]. It is hence of interest to see that for $\kappa \leq 1$, $\hat{\lambda}/\lambda \to_p 1$ holds whenever the likelihood ratio test succeeds (asymptotically) in rejecting the global null hypothesis.

**4. Numerical examples.** It emerged from the analysis above that the standard deviation–proportional bounding function is optimal in an asymptotic sense. In the following discussion we briefly compare various bounding functions for a moderate number of tests, $n = 1000$. The setup is identical to the shift-location testing of Section 3.2, equation (18). For true null hypotheses, test statistics follow the normal distribution $\mathcal{N}(0, 1)$. For false null hypotheses, test statistics are shifted by an amount $\mu > 0$ and are $\mathcal{N}(\mu, 1)$-distributed.

The proportion $\hat{\lambda}/\lambda$ of correctly identified false null hypotheses is computed for various values of the shift parameter $\mu$ and three bounding functions. The results for 100 simulations are shown in Figure 3. The left column shows results for very few false null hypotheses ($\lambda = 0.01$), corresponding to 10 false null hypotheses, while results are shown in the right column for a moderately large number of false null hypotheses ($\lambda = 0.2$).

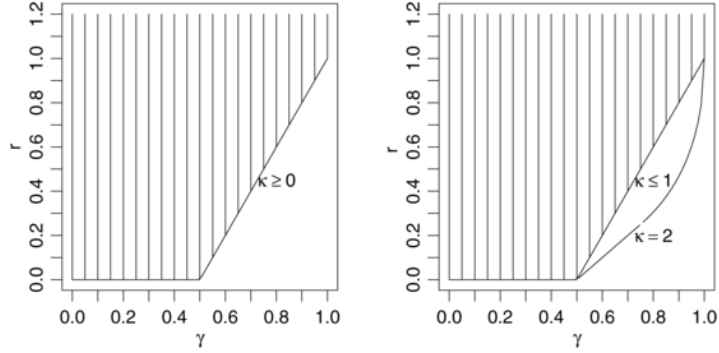

FIG. 2. *Comparison between the estimate of $\hat{\lambda}$ and detection regions under higher criticism if test statistics follow the location-shift model* (18) *and are distributed according to the generalized Gaussian Subbotin distribution with shift parameter* (19). *The shaded area in the left panel shows again the area of convergence in probability of $\hat{\lambda}/\lambda$ to 1 for a bounding function in the class $Q_{1/2}$. The shaded area in the right panel corresponds to the region where higher criticism can reject asymptotically the null hypothesis $H_0 : \overline{\lambda} = 0$ for $\kappa \leq 1$, including the double-exponential case. The line below marks the detection boundary for the Gaussian case ($\kappa = 2$).*

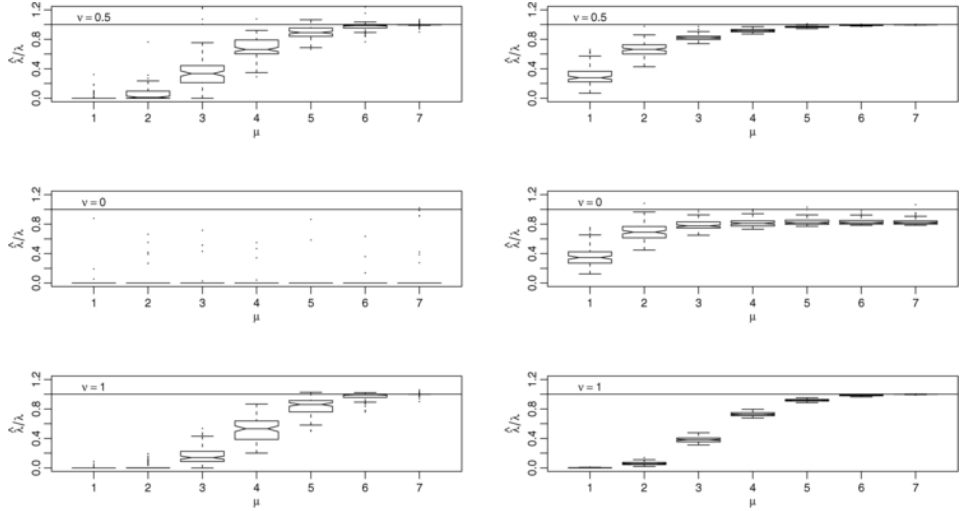

FIG. 3. *The proportion $\hat{\lambda}/\lambda$ of correctly detected false null hypotheses as a function of the separation $\mu$. Results are shown for the standard deviation–proportional bounding function (top row), the constant bounding function (middle row), and the linear bounding function (bottom row).*

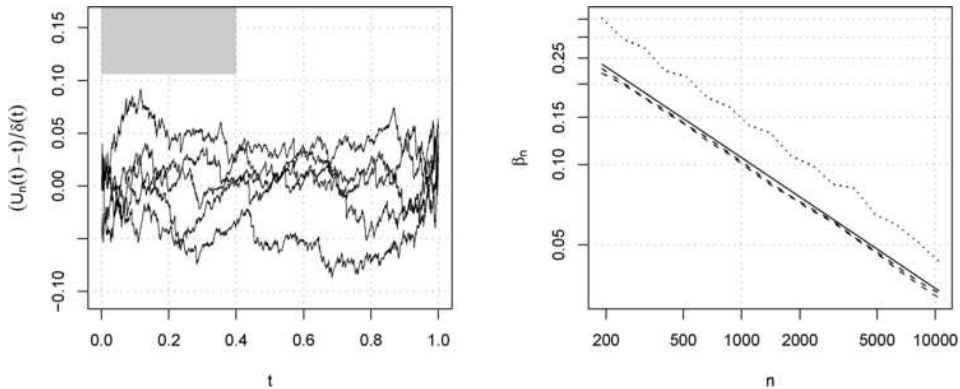

Fig. 4. *Random samples of the weighted empirical distribution function $(U_n(t) - t)/\delta(t)$ with $\delta(t) = \sqrt{t(1-t)}$ on the left. Various bounding sequences $\beta_{n,\alpha}$ as a function of $n$ in log-log scale on the right: the asymptotically valid bounding sequence (solid line), and the bounding sequences for the intervals $(0, 1)$ (dotted line), $(1/n, 1 - 1/n)$ (upper dashed line) and $(1/n, 0.01)$ (lower dashed line), as obtained by simulation. Note that the latter two are almost indistinguishable.*

For very few false null hypotheses ($\lambda = 0.01$), both the standard deviation–proportional and linear bounding functions identify a substantial proportion of false null hypotheses if the shift $\mu$ is larger than about 3. The expected value of the largest test statistic from true null hypotheses is, for comparison in the current setup, at around 3.7. The constant bounding function ($\nu = 0$) fails to identify any of the 10 false null hypotheses even for very large shifts $\mu$. This is in line with the theoretical results from Section 3.2. For a moderately large number of false null hypotheses ($\lambda = 0.2$), the performance of the linear bounding function is worse than for the other two bounding functions, as expected from the asymptotic results in Section 3.1. The standard deviation–proportional bounding function ($\nu = 1/2$) in both cases consistently identifies the most false null hypotheses, and the optimality of this bounding function is thus numerically evident for moderate sample sizes as well.

For the standard deviation–proportional bounding function ($\nu = 1/2$), asymptotic control was proposed in (10). The result relies on convergence of the supremum of a weighted empirical distribution to the Gumbel distribution. This convergence is in general slow, as already mentioned in Remark 3. The convergence is comparably fast, however, if the region over which the supremum is taken is restricted to, say, $(1/n, 1 - 1/n)$, as observed by Donoho and Jin [6]. We illustrate this in the following text. Restricting the interval over which the supremum is taken in (5) to some interval $(a, b)$ with $0 < a < b < 1$, bounding sequences can be defined analogous to Definition 1 by requirement (b) in Definition 1 and

(25) $$\beta_{n,\alpha} = \min\left\{\beta : P\left(\sup_{t\in(a,b)} \frac{U_n(t) - t}{\delta(t)} > \beta\right) \leq \alpha\right\}.$$

Bounding sequences for the interval $(0,1)$ satisfy (25) for every interval $(a,b)$, but might be unduly conservative. Less conservative bounding sequences can be found conveniently by approximating the probability of $\sup_{t\in(a,b)}(U_n(t)-t)/\delta(t) > \beta$ with the empirical proportion of occurrence of this event among a large number of simulations. This is illustrated in the left panel in Figure 4. Shown are five random samples of the the weighted empirical distribution $(U_n(t)-t)/\delta(t)$ for $n=200$ and $\delta(t)=\sqrt{t(1-t)}$. Let the value $\beta$ correspond to the lower bound of the gray area in Figure 4. For an interval $(a,b)=(0,0.4)$, the event $\sup_{t\in(a,b)}(U_n(t)-t)/\delta(t) > \beta$ corresponds then to the event that a realization of a weighted empirical distribution crosses the gray area. The bounding sequences obtained by using 1000 simulations of the weighted empirical distribution are shown in the right panel in Figure 4 for various intervals $(a,b)$.

There are two main conclusions. First, one might suspect that $p$-values from false null hypotheses are mostly found in a neighborhood around zero. Restricting the region in (5) to such a neighborhood promises thus to capture all $p$-values from false null hypotheses while allowing for smaller bounding sequences. However, the numerical results suggest otherwise. The bounding sequence for the region $(1/n, 0.01)$ is, for example, almost indistinguishable from the bounding sequence for the region $(1/n, 1-1/n)$, as can be seen in Figure 4.

Second, the agreement of the asymptotically valid bounding sequence (13) with the bounding sequence that is obtained by simulation for the interval $(1/n, 1-1/n)$ is very good even for moderate sample sizes, while the agreement is not so good for the interval $(0,1)$. When using the asymptotically valid bounding sequence it is hence advisable to restrict the region over which the supremum is taken in (5) to $(1/n, 1-1/n)$. This ensures that the true level is close to the chosen level $\alpha$ for moderate sample sizes.

For practical applications, we hence recommend that one calculate the supremum in (5) over a region $(1/n, 1-1/n)$ and use the standard deviation–proportional bounding function with the asymptotically valid bounding sequence (13). The asymptotic results of the previous sections hold for this modified procedure.

## 5. Proofs.

PROOF OF THEOREM 2. First it is shown that, as long as $\gamma \in (0, \frac{1}{2})$ and $\nu \leq \frac{1}{2}$, for any given $\varepsilon > 0$,

(26) $$P(\hat{\lambda} < (1-\varepsilon)\lambda) \to 0, \qquad n \to \infty.$$

Let the empirical distribution of $p$-values be defined as in (4) by $F_n(t) = n^{-1} \sum_{i=1}^{n} \mathbf{1}\{P_i \leq t\}$. We suppose that the proportion of false null hypotheses is fixed at $\lambda$, so that $F_n(t)$ is a mixture $F_n(t) = \lambda G_{n_1}(t) + (1-\lambda)U_{n_0}(t)$, where $G_{n_1}(t)$ is the empirical distribution of $n_1 = \lambda n$ i.i.d. $p$-values with distribution $G$ and $U_{n_0}(t)$ is the empirical distribution of $n_0 = (1-\lambda)n$ i.i.d. $p$-values with uniform distribution $U$. For any $t < 1$,

$$\hat{\lambda} = \sup_{t \in (0,1)} \frac{F_n(t) - t - \beta_{n,\alpha_n}\delta(t)}{1-t} \tag{27}$$

$$\geq \frac{F(t) - t}{1-t} + \frac{F_n(t) - F(t) - \beta_{n,\alpha_n}\delta(t)}{1-t} \tag{28}$$

$$= \lambda \frac{G(t) - t}{1-t} + \frac{F_n(t) - F(t) - \beta_{n,\alpha_n}\delta(t)}{1-t}. \tag{29}$$

Whereas $\inf_{t \in (0,1)} G'(t) = 0$ and, hence, $\sup_{t \in (0,1)} (G(t) - t)/(1-t) = 1$, there exists by continuity of $G(t)$ some $t_1$ so that $(G(t_1) - t_1)/(1-t_1) > (1-\varepsilon/2)$. Setting $\tilde{\varepsilon} = \frac{1}{2}(1-t_1)\varepsilon$, it suffices to show that for every $\varepsilon > 0$,

$$P(\beta_{n,\alpha_n}\delta(t_1) + F(t_1) - F_n(t_1) > \tilde{\varepsilon}\lambda) \to 0, \qquad n \to \infty.$$

Whereas $F_n(t_1) - F(t_1) = O_P(n^{-1/2})$ and $\lambda \sim n^{-\gamma}$ with $\gamma < \frac{1}{2}$, this follows from the finiteness of $\delta(t)$ and, because $\alpha_n$ vanishes slowly, from (17). This completes the first part of the proof of Theorem 2.

For the second part, it suffices to show that for $\nu \in (\frac{1}{2}, 1]$ and $\gamma \in (1-\nu, \frac{1}{2})$, and any $\varepsilon > 0$,

$$P(\hat{\lambda} > \varepsilon\lambda) \to 0, \qquad n \to \infty. \tag{30}$$

In this regime, the penalty $\beta_{n,\alpha_n}\delta(t)$ is asymptotically larger than the signal from false null hypotheses. Using the definition of $\hat{\lambda}$, the notation $n_0 = (1-\lambda)n$ and $n_1 = \lambda n$, and $F_n(t) = \lambda G_{n_1}(t) + (1-\lambda)U_{n_0}(t)$, it follows that

$$P(\hat{\lambda} > \varepsilon\lambda) = P\left(\sup_{t \in (0,1)} \frac{F_n(t) - t - \beta_{n,\alpha_n}\delta(t)}{1-t} > \varepsilon\lambda\right)$$

$$= P\left(\sup_{t \in (0,1)} \lambda \frac{G_{n_1}(t) - t}{1-t} - \varepsilon\lambda + (1-\lambda)\frac{U_{n_0}(t) - t}{1-t} - \beta_{n,\alpha_n}\frac{\delta(t)}{1-t} > 0\right)$$

$$\leq P\left(\sup_{t \in (0,1)} \lambda \frac{G_{n_1}(t) - t}{1-t} - \varepsilon\lambda - \frac{\beta_{n,\alpha_n}}{2}\frac{\delta(t)}{1-t} > 0\right) \tag{31}$$

$$+ P\left(\sup_{t \in (0,1)} (1-\lambda)\frac{U_{n_0}(t) - t}{1-t} - \frac{\beta_{n,\alpha_n}}{2}\frac{\delta(t)}{1-t} > 0\right). \tag{32}$$

Observe in (32) that $(1-\lambda)^{-1}\beta_{n,\alpha_n} = n\beta_{n,\alpha_n}/n_0 \geq \beta_{n_0,\alpha_n} \geq \beta_{n_0,\alpha_{n_0}}$. Thus (32) can be bounded by $P(V_{n_0,\delta} > \beta_{n_0,\alpha_{n_0}}/2)$. By (16) and $n_0 \to \infty$ it follows that (32) vanishes for $n \to \infty$. It remains to show that (31) vanishes as well. Let $t_2 = \sup\{t \in (0,1) : G(t) \le \varepsilon/2\}$. Using Bonferroni's inequality, (31) is bounded by

$$P\left(\sup_{t \in (0,t_2]} \lambda \frac{G_{n_1}(t) - t}{1-t} - \varepsilon\lambda > 0\right) \tag{33}$$

$$+ P\left(\sup_{t \in (t_2,1)} \lambda \frac{G_{n_1}(t) - t}{1-t} - \frac{\beta_{n,\alpha_n}}{2} \frac{\delta(t)}{1-t} > 0\right). \tag{34}$$

Whereas $(G_{n_1}(t) - t)/(1-t) \le G_{n_1}(t)$ for all $t \in [0,1]$, the first term (33) is bounded by $P(G_{n_1}(t_2) > \varepsilon)$, which vanishes for $n \to \infty$ because by definition of $t_2$, $G(t_2) \le \varepsilon/2$ and $n_1 = \lambda n \to \infty$. Using $G_{n_1}(t) \le 1$, the second term (34) equals zero if $\beta_{n,\alpha_n} \inf_{t \in (t_2,1)} \delta(t)/(1-t) > 2\lambda$. By conditions (a) and (b) in Definition 3, it holds that $\inf_{t \in (t_2,1)} \delta(t)/(1-t) > 0$. By (14), it follows furthermore that $\beta_{n,\alpha_n}/\lambda \to \infty$ for $n \to \infty$, which completes the proof. □

PROOF OF THEOREM 3. First it is shown that for $r > \frac{1}{\nu}(\gamma - \frac{1}{2})$,

$$P(\hat{\lambda} < (1-\varepsilon)\lambda) \to 0, \qquad n \to \infty. \tag{35}$$

Here the penalty is again asymptotically larger than the signal from false null hypotheses for a fixed point $t \in (0,1)$. However, because the signal from false null hypotheses is increasing in strength for larger $n$, the evidence for a certain amount of false null hypotheses can be found at decreasing values of $t$. Using the definition of $\hat{\lambda}$, for any $t \in (0,1)$, $\hat{\lambda} \ge F_n(t) - t - \beta_{n,\alpha_n} \delta(t)$ and, hence, for any $t \in (0,1)$,

$$\hat{\lambda}/\lambda - 1 \ge (1 - G_{n_1}(t)) - t - \frac{1-\lambda}{\lambda}(t - U_{n_0}(t)) - \frac{1}{\lambda}\beta_{n,\alpha_n}\delta(t),$$

where again $n_1 = \lambda n$ and $n_0 = (1-\lambda)n$. Choosing $t_{n,\tau} = n^{-r+\tau}$ for some $0 < \tau < r - \frac{1}{\nu}(\gamma - \frac{1}{2})$, observe that by (20) it follows that $1 - G(n^{-r+\tau}) = o(1)$. Hence

$$\hat{\lambda}/\lambda - 1 \ge (1 - G(t_{n,\tau})) - |G(t_{n,\tau}) - G_{n_1}(t_{n,\tau})|$$
$$- t_{n,\tau} - \lambda^{-1}|t_{n,\tau} - U_{n_0}(t_{n,\tau})| - \lambda^{-1}\beta_{n,\alpha_n}\delta(t_{n,\tau})$$
$$= o(1) - o_p(1) - o(1) - O_p(n^{\gamma - (1/2 + (r-\tau)/2)})$$
$$- O(n^{\gamma - (1/2 + (r-\tau)\nu)} \log n).$$

The proof of (35) follows because $\gamma < \frac{1}{2} + \nu(r - \tau) \le \frac{1}{2} + \frac{r-\tau}{2}$.

Second, it has to be shown that $P(\hat{\lambda} > \varepsilon\lambda) \to 0$ if $r < \frac{1}{\nu}(\gamma - \frac{1}{2})$. Again, the evidence for a certain amount of false null hypotheses would have to be found at decreasing values of $t$. However, the decrease has to be so fast in this regime that the signal from false null hypotheses is not captured. Using again the notation $n_1 = \lambda n$ and $n_0 = (1-\lambda)n$, we find that $\hat{\lambda} = \sup_{t \in (0,1)} D_{n,\lambda}(t)$, where

$$(36) \quad D_{n,\lambda}(t) := \frac{\lambda(G_{n_1}(t) - t) + (1-\lambda)(U_{n_0}(t) - t) - \beta_{n,\alpha_n}\delta(t)}{1 - t}.$$

Choose a sequence $t_{n,\rho} = n^{-r-\rho}$ for some $0 < \rho < \frac{1}{\nu}(\gamma - \frac{1}{2}) - r$. The regions $(0, t_{n,\rho}]$ and $(t_{n,\rho}, 1)$ are considered separately for the following. In particular, it is shown that both $P(\sup_{t \in (0,t_{n,\rho}]} D_{n,\lambda}(t) > \varepsilon\lambda)$ and $P(\sup_{t \in (t_{n,\rho},1)} D_{n,\lambda}(t) > \varepsilon\lambda)$ vanish for $n \to \infty$. For $t > t_{n,\rho}$, it holds that

$$P\bigg(\sup_{t \in (t_{n,\rho},1)} D_{n,\lambda}(t) > \varepsilon\lambda\bigg)$$

$$\leq P\bigg(\sup_{t \in (t_{n,\rho},1)} \lambda + (1-\lambda)\frac{U_{n_0}(t) - t}{1 - t} - \beta_{n,\alpha_n}\frac{\delta(t)}{1 - t} > 0\bigg)$$

$$\leq P\bigg(\sup_{t \in (t_{n,\rho},1)} (1-\lambda)\frac{U_{n_0}(t) - t}{1 - t} - \frac{\beta_{n,\alpha_n}}{2}\frac{\delta(t)}{1 - t} > 0\bigg)$$

$$+ \mathbf{1}\bigg\{\sup_{t \in (t_{n,\rho},1)} \lambda - \frac{\beta_{n,\alpha_n}}{2}\frac{\delta(t)}{1 - t} > 0\bigg\}$$

$$(37) \quad = P\bigg(\sup_{t \in (t_{n,\rho},1)} (U_{n_0}(t) - t) - \frac{n}{n_0}\frac{\beta_{n,\alpha_n}}{2}\delta(t) > 0\bigg)$$

$$(38) \quad + \mathbf{1}\bigg\{\inf_{t \in (t_{n,\rho},1)} \frac{\beta_{n,\alpha_n}}{2}\frac{\delta(t)}{1 - t} < \lambda\bigg\}.$$

By (16) and because $n\beta_{n,\alpha_n}$ is monotonically increasing, (37) vanishes for $n \to \infty$. For (38), because $\delta \in Q_\nu$, there exists a constant $c$ so that $\inf_{t \in (t_{n,\rho},1)} \delta(t_{n,\rho}) \geq cn^{-\nu(r+\rho)}$. It follows by $r + \rho < \frac{1}{\nu}(\gamma - \frac{1}{2})$ that $\inf_{t \in (t_{n,\rho},1)} \beta_{n,\alpha_n}\delta(t_{n,\rho})/\lambda \to \infty$ for $n \to \infty$, which completes the first part of the proof.

It remains to show that $P(\sup_{t \in (0,t_{n,\rho}]} D_{n,\lambda}(t) > \varepsilon\lambda) \to 0$ for $n \to \infty$. It holds that

$$P\bigg(\sup_{t \in (0,t_{n,\rho}]} D_{n,\lambda}(t) > \varepsilon\lambda\bigg)$$

$$(39) \quad \leq P\bigg(\sup_{t \in (0,t_{n,\rho}]} (1-\lambda)\frac{U_{n_0}(t) - t}{1 - t} - \beta_{n,\alpha_n}\frac{\delta(t)}{1 - t} > \frac{\varepsilon}{3}\lambda\bigg)$$

$$(40) \quad + P\bigg(\sup_{t \in (0,t_{n,\rho}]} \lambda\frac{G_{n_0}(t) - G(t)}{1 - t} > \frac{\varepsilon}{3}\lambda\bigg)$$

$$(41) \quad + \mathbf{1}\bigg\{\sup_{t \in (0,t_{n,\rho}]} \lambda\frac{G(t) - t}{1 - t} > \frac{\varepsilon}{3}\lambda\bigg\}.$$

As already argued above, the probability on the right-hand side of (39) vanishes for $n \to \infty$. The probability (40) clearly likewise vanishes and it remains to show that (41) vanishes as well for $n \to \infty$. Whereas $t_{n,\rho} \to 0$, it holds that $(1-t)^{-1} \leq 2$ for $t \in (0, t_{n,\rho}]$ and large enough values of $n$. The term (41) vanishes hence if $G(t_{n,\rho}) < \frac{\varepsilon}{6}$. This is equivalent to $\log G^{-1}(\frac{\varepsilon}{6}) < -(\tau + \rho) \log n$, and the claim follows from property (20). $\square$

**Acknowledgments.** The authors would like to thank an anonymous referee, the Associate Editor and Jianqing Fan for helpful comments, which helped to improve an earlier version of the paper. Nicolai Meinshausen would also like to thank Peter Bühlmann for interesting discussions.

Seminar für Statistik  
ETH-Zürich  
8092 Zürich  
Switzerland  
E-mail: nicolai@stat.math.ethz.ch

Department of Statistics  
University of California  
Berkeley, California 94720  
USA  
E-mail: rice@stat.berkeley.edu